\documentclass{amsart}
\usepackage{amssymb, amsthm, amsmath}
\newtheorem{thm}{Theorem}[section]
\newtheorem{prop}{Proposition}[section]

\newtheorem{cor}{Corollary}[section]
\theoremstyle{definition}
\newtheorem{prob}{Problem}[section]
\newtheorem{dfn}{Definition}[section]
\newtheorem{rem}{Remark}[section]
\newtheorem{ex}{Example}[section]
\def\m{{\mathfrak{m}}}
\def\Hom{{\rm{Hom}}}
\def\Soc{{\rm{Soc}}}
\begin{document}
\title[Lefschetz elements of Gorenstein algebras]{Lefschetz elements of Artinian Gorenstein algebras and Hessians 
of homogeneous polynomials}
\author[T. Maeno]{Toshiaki Maeno}
\address{Toshiaki Maeno, Department of Electrical Engineering, 
Kyoto University 
Kyoto 606-8501, Japan}
\email{maeno@kuee.kyoto-u.ac.jp}
\thanks{The first author is supported by Grant-in-Aid for Scientific Research.}
\author[J. Watanabe]{Junzo Watanabe}
\address{Junzo Watanabe, Department of Mathematics, 
Tokai University 
Hiratsuka 259-1292, Japan}
\email{junzowat@keyaki.cc.u-tokai.ac.jp}
\date{}
\begin{abstract}
We give a characterization of the Lefschetz elements 
in Artinian Gorenstein rings over a field of characteristic zero 
in terms of the higher Hessians. As an application, we give new 
examples of Artinian Gorenstein rings which do not have 
the strong Lefschetz property. 
\end{abstract}
\maketitle
\section*{Introduction}
The Lefschetz property is a ring-theoretic abstraction of the Hard Lefschetz 
Theorem for compact K\"ahler manifolds (see e.g. \cite{GH}). 
The following are fundamental problems on the study 
of the Lefschetz property for Artinian graded algebras: 
\begin{prob}
For a given graded Artinian algebra $A,$ decide whether or not $A$ has the 
strong (or weak) Lefschetz property.
\end{prob}
\begin{prob}
When a graded Artinian algebra $A$ has the strong Lefschetz property, 
determine the set of Lefschetz elements in $A_1.$ 
\end{prob}
In \cite{W1}, it was shown that ``most'' Artinian Gorenstein 
algebras have the strong Lefschetz property. However, it is a difficult 
problem to know whether a given graded Artinian algebra has the strong (or weak) 
Lefschetz property. 
In principle, if a graded Artinian algebra $A$ over a field $k$ is given with a presentation 
\[ A=k[x_1,\ldots,x_n]/(f_1,\ldots,f_m), \] 
we have an algorithm to answer the 
above Problems since it is sufficient to compute the determinants of 
the matrix expression of the multiplication map by a general element of $A_1$ 
with respect to an arbitrary homogeneous linear basis of $A.$ In particular, 
the complement of the set of Lefschetz elements in $A_1$ 
has a structure of the algebraic set defined by certain determinants. 
However, it is hard in general to carry out the computation based on this algorithm 
even with the help of computer. 

In the present paper, we give a simple criterion to answer these problems 
for Artinian Gorenstein algebras over a field $k$ of characteristic zero. 
It is known that a graded Artinian Gorenstein algebra is characterized by the ``Poincar\'e 
duality'' which holds for the cohomology ring of the compact oriented 
manifolds. 
Hence graded Artinian Gorenstein algebras with the strong Lefschetz property 
are a natural class of commutative 
algebras comparable to the cohomology ring of compact K\"ahler manifolds. 

A typical example of graded Artinian Gorenstein algebras is the coinvariant 
algebra of finite Coxeter groups. In fact, the coinvariant algebra of the Weyl 
group is isomorphic to the cohomology ring of the corresponding flag variety. 
In \cite{MNW} and \cite{NW}, it has been shown that the coinvariant algebra 
of any finite Coxter group has the Lefschetz property and that the set of the Lefschetz 
elements is the complement of the union of the reflection hyperplanes except for 
type $H_4$ case. The determination of the set of the Lefschetz elements 
is still open for $H_4$ because of the computational complexity.  

Let us consider the polynomial ring $R=k[x_1,\ldots,x_n]$ and 
the algebra of differential operators 
\[ Q=k[ \frac{\partial}{\partial x_1},\ldots,\frac{\partial}{\partial x_n} ]. \]  
Every graded Artinian Gorenstein algebra has the presentation 
\[ A\cong Q/{\rm Ann}_Q F, \quad 
{\rm Ann}_QF= \{ \varphi(\partial_1,\ldots,\partial_n) \in Q \; | \; 
\varphi(\partial_1,\ldots,\partial_n)F(x)=0 \} , \] 
for some homogeneous polynomial $F \in k[x_1,\ldots,x_n].$ 
We introduce the higher Hessians ${\rm Hess}^{(d)}F,$ $1\leq d \leq [\deg F /2],$ 
of the polynomial $F$ 
in order to describe the condition for an element $L \in A_1$ to be a strong 
Lefschetz element. 
The set of the strong Lefschetz elements of $A$ is a Zariski open set in $A_1,$ 
which is given as the complement of all the zero loci of the higher Hessians. 
We will discuss the explicit description of the set of Lefschetz elements 
of $A=Q/{\rm Ann}_Q F$ for the Fermat type polynomial 
$F=\sum_{i=1}^nx_i^n-n(n-1)s\prod_{i=1}^nx_i.$ 

When one of the higher Hessians of $F$ is identically zero, the algebra 
$Q/{\rm Ann}_QF$ does not have the strong Lefschetz property. In \cite{HMNW}, \cite{Ik} and 
\cite{W2}, examples of Artinian Gorenstein algebras which do not have the strong 
Lefschetz property are given. The examples in \cite{HMNW} and \cite{W2} are based on the polynomials 
with the zero Hessian. 
In the last section, we give some polynomials $F$ such that 
${\rm Hess}\; F\not=0$ and ${\rm Hess}^{(2)}F=0$ to get new examples 
of Artinian Gorenstein algebras which do not have the strong Lefschetz property. 
\subsection*{Acknowledgements} The first author would like to thank Akihito Wachi 
for informing of Proposition 2.1 and for useful comments. 
The authors are grateful to Tadahito Harima and Anthony Iarrobino for their suggestions and comments. 

\section{Lefschetz properties}
\begin{dfn}
Let $A=\oplus_{d=0}^DA_d,$ $A_D\not=0,$ be a graded Artinian algebra. \\ 
(1) We say that $A$ has the strong Lefschetz property if there exits an element 
$L \in A_1$ such that the multiplication map 
\[ \times L^d:A_i \rightarrow A_{i+d} \]  
is of full rank $($i.e. injective or surjective$)$ for all $0\leq i\leq D$ and 
$0\leq d\leq D-i.$ 
We call $L\in A_1$ with this property a strong Lefschetz element. \\ 
(2) If we assume the existence of $L\in A_1$ 
such that 
\[ \times L:A_i \rightarrow A_{i+1} \] 
is of full rank for $i=0,\ldots,D-1,$ we say 
that $A$ has the weak Lefschetz property. 
\end{dfn}

If a graded Artinian algebra $A$ over a field $k$ is generated by $A_1$ as a $k$-algebra, 
we say that $A$ has the standard grading. The weak Lefschetz 
property implies the unimodality of the Hilbert function, provided that 
the $k$-algebra $A$ has the standard grading. 

\begin{dfn}
Let $A=\oplus_{d=0}^DA_d,$ $A_D\not=0,$ be a graded Artinian algebra. 
We say that $A$ has the strong Lefschetz property in the narrow sense if there exits an element 
$L \in A_1$ such that the multiplication map 
\[ \times L^{D-2i}:A_i \rightarrow A_{D-i} \] 
is bijective for $i=0,\ldots,[D/2].$ 
\end{dfn}

If a graded Artinian $k$-algebra $A$ has the strong Lefschetz property 
in the narrow sense, then the Hilbert function of $A$ is unimodal and symmetric. 
When a graded Artinian $k$-algebra $A$ has a symmetric Hilbert function, 
the notion of the strong Lefschetz property on $A$ coincides with the one in the 
narrow sense.  Our main interest in this paper is to consider Artinian 
{\it Gorenstein} algebras, so the strong Lefschetz property will be used 
in the narrow sense in the subsequent sections. 
Throughout this paper, graded Artinian $k$-algebras $A=\oplus_{d=0}^DA_d$ are assumed 
to satisfy the conditions $A_0 \cong k$ and $\dim_k A_1, \dim_k A_D >0.$ 

\section{Artinian Gorenstein algebra} 
Throughout $k$ denotes a field. 
\begin{dfn} (See \cite[Chapter 5, 6.5]{Sm}.) 
A finite-dimensional graded $k$-algebra $A=\oplus_{d=0}^DA_d$ is called the Poincar\'e 
duality algebra if $\dim_k A_D=1$ and the bilinear pairing 
\[ A_d \times A_{D-d} \rightarrow A_D \cong k \] 
is non-degenerate for $d=0,\ldots,[D/2].$ 
\end{dfn}
The following is a well-known fact (see e.g. \cite{GHMS}).
\begin{prop}
A graded Artinian $k$-algebra $A$ is a Poincar\'e duality algebra if and only if 
$A$ is Gorenstein. 
\end{prop}
\begin{proof} Assume that $A=\oplus_{d=0}^DA_d$ is a Poincar\'e duality algebra. 
For any element $f\in A\setminus \{ 0 \}$ of degree less than $D,$ 
there exisits an element 
$g\in A\setminus A_0$ such that $fg\not=0.$ 
Hence the socle ideal $\Soc(A)$ of $A$ coincides with 
a one-dimensional $k$-subspace $A_D.$ This means that $A$ is Gorenstein. Conversely, 
if $A$ is Gorenstein, the socle ideal $\Soc(A)$ is the one-dimensional $k$-vector space. 
Since the maximal degree part $A_D$ is contained in $\Soc(A),$ 
we have that $A_D=\Soc(A)$ and $\dim_kA_D=1.$ 
We will prove the following claim by induction on $d$: \smallskip \\ 
$(*)_d$ \hspace{1cm} if $f\in A_{D-d}$ satisfies $fg=0$ for all $g\in A_d,$ then $f=0.$ \smallskip \\ 
If $f\in A_{D-1}$ satisfies $fg=0$ for all $g\in A_1,$ 
then $fg=0$ for all $g\in A_{>0}$ for a degree reason. 
This implies that $f\in A_{D-1} \cap \Soc(A)=0,$ 
so $(*)_1$ follows. 
Let us assume that a nonzero element $f\in A_{D-d}\setminus\{ 0 \},$ $d>1,$ satisfies $fg=0$ for 
all $g\in A_d,$ and that there exists an element $h\in A_i,$ $1\leq i <d,$ 
such that $\varphi:=fh\not=0.$ 
By the induction hypothesis $(*)_{d-i},$ 
we can find 
an element $h'\in A_{d-i}$ such that $\varphi h' \not=0$ for the nonzero 
element $\varphi \in A_{D-d+i}.$ 
Then we have $\varphi h' =f(hh')\not=0,$ which  
is a contradiction since $hh'\in A_d.$  
We have proved that if $f\in A_{D-d}$ satisfies 
$fg=0$ for all $g\in A_d,$ then we have $f=0$ by contradiction. 
Now the claim $(*)_d$ is proved. 
The claims $(*)_d$ for $d=1,\ldots,D$ imply that 
the pairing $A_d \times A_{D-d} \rightarrow A_D$ is non-degenerate for $d=1,\ldots,[D/2].$  
\end{proof}
\begin{rem} 
(1) The above proposition shows that the even part of the cohomology ring $H^{even}(M,k)$ 
with coefficient in a field $k$ of characteristic zero 
of any compact orientable manifold $M$ is Gorenstein. \\ 
(2) The Poincar\'e duality algebra is an abstraction of the property 
of the cohomology ring of compact orientable manifolds, whereas the strong Lefschetz 
property is inspired by the Hard Lefschetz Theorem for compact K\"ahler 
manifolds. Though the K\"ahler manifold is always oriented, the strong Lefschetz 
property does not imply the Poincar\'e duality. In other words, there exist examples of 
graded Artinian {\it non}-Gorenstein algebras with the strong Lefschetz property. 
For example, $A=k[x,y]/(x^2,xy,y^3)$ is a non-Gorenstein algebra with the strong 
Lefschetz property. At the same time, the Poincar\'e duality does not imply the strong 
Lefschetz property. See Example 2.1 and Section 5. 
\end{rem}

For simplicity, we assume that the characteristic of the field $k$ to be zero 
in the rest of this paper, though our main results hold also when 
the characteristic of $k$ is greater than the socle degree $D$ of the Gorenstein 
algebra $A.$ 
Let us regard the polynomial algebra $R:=k[x_1,\ldots,x_n]$ as a module 
over the algebra $Q:= k[X_1,\ldots,X_n]$ via the identification $X_i=\partial / \partial x_i.$ 
For a polynomial $F\in R,$ we define the ideal ${\rm Ann}_Q F$ of $Q$ by 
\[ {\rm Ann}_Q F:= \{ a(X_1,\ldots,X_n)\in Q \; | \; a(\partial_1,\ldots,\partial_n)F=0 \}. \] 

The following theorem is a well-established fact among the experts. 
In fact it is immediate if the theory of the inverse system is 
taken for granted (\cite{BH},\cite{Ge},\cite{GW}). 
However the theory of inverse system does not seem to be well-known to non-specialists, 
so we give a direct proof for it. 
\begin{thm} 
Let $I$ be an ideal of $Q=k[X_1,\ldots,X_n]$ and 
$A=Q/I$ the quotient algebra. Denote by $\m$ the maximal ideal 
$(X_1,\ldots,X_n)$ of $Q.$ 
Then $\sqrt{I}=\m$ and the $k$-algebra $A$ is Gorenstein if and only if 
there exists a polynomial $F \in R=k[x_1,\ldots,x_n]$ 
such that $I={\rm Ann}_Q F.$ 
\end{thm}
\begin{proof} Assume that $I={\rm Ann}_Q F$ for some polynomial 
$F \in R.$ Since $F$ is annihilated by differential operators of sufficiently 
high order, ${\rm Ann}_Q F$ contains $\m^l$ for sufficiently large $l.$ 
Since we are working over a field $k$ of characteristic zero, 
it is clear that there exists a polynomial $G\in Q$ such that 
$G(X)F(x) \in k^{\times}.$ We will show in the following that $I:\m = I+k\cdot G.$ 
Since $G(X)F(x)$ is a constant, it immediately follows that 
\[ \partial_i G(\partial_1,\ldots,\partial_n)F(x_1,\ldots,x_n)=0, \quad 
i=1,\ldots,n. \] 
This shows that $G\in I:\m.$ Now let $a(X_1,\ldots,X_n)\in I:\m$ be any element. 
By definition of the ideal $I:\m,$ we have $X_ia(X_1,\ldots,X_n)\in I={\rm Ann}_Q F$ for 
all $i=1,\ldots,n.$ 
This means that 
\[ \partial_i a(\partial_1,\ldots,\partial_n)F(x_1,\ldots,x_n)=0, \quad 
i=1,\ldots,n. \] 
Hence we have that $a(X)F(x)$ is a constant. As we have already seen that 
$G(X)F(x)$ is a non-zero constant, we have that $a(X)-cG(X) \in {\rm Ann}_Q F$ 
for some constant $c\in k.$ We have shown that $I:\m = I+k\cdot G.$ In other words, 
the $k$-vector space $(I:\m)/I$ is one-dimensional. Thus $A=Q/{\rm Ann}_Q F$ is 
a Gorenstein algebra (see e.g. \cite{BH}). 

Now let us prove the converse implication. Assume that $A=Q/I$ is an Artinian 
Gorenstein algebra. Then we have the isomorphism $\Hom_k(A,k) \cong A$ as an $A$-module. 
The $Q$-module $\Hom_k(Q,k)$ is identified with the ring of formal power series 
$\hat{R}=k[[ x_1,\ldots,x_n]]$ regarded as a $Q$-module. From the exact sequence 
$Q \rightarrow A \rightarrow 0,$ we have the exact sequence of $Q$-modules: 
\[ 0 \rightarrow \Hom_k(A,k) \cong A \stackrel{\theta}{\rightarrow} \Hom_k(Q,k) \cong 
k[[ x_1,\ldots,x_n]]. \] 
Define $F \in k[[ x_1,\ldots,x_n]]$ as the image of $1\in A$ 
by the homomorphism $\theta.$ From the assumption that $I$ contains $\m^l$ for $l\gg 0,$ 
the image of $\theta$ annihilates polynomials in $Q$ of sufficiently large 
degrees, so $F$ is a polynomial in $R.$ Finally we have that 
\[ {\rm Ann}_Q F = \{ a\in Q \; | \; ab\in I, \; \forall b \in Q \} = I, \] 
so $A=Q/{\rm Ann}_Q F.$ 
\end{proof}
\begin{rem} 
Let $k$ be the field ${\bf C}$ of complex numbers. In this case, 
for a polynomial 
\[ F=\sum c_{i_1\ldots i_n}x_1^{i_1}\cdots x_n^{i_n}, \quad c_{i_1\ldots i_n}\in {\bf C}, \] 
we can choose the complex conjugate 
\[ \overline{F}=\sum \overline{c}_{i_1\ldots i_n}X_1^{i_1}\cdots X_n^{i_n} \] 
as a generator of the socle of $A=Q/{\rm Ann}_Q F.$ 
\end{rem}
\begin{rem}
When $I$ is a homogeneous ideal (i.e. $A$ is graded), 
the condition $\sqrt{I}=\m$ is satisfied. 
In this case, 
we can choose $F$ as a homogeneous polynomial. 
\end{rem}
\begin{ex}
Stanley \cite{S2} gave an example of Artinian Gorenstein algebra 
with a non-unimodal Hilbert function. 
Let us take the polynomial 
\[ F(u,v,w,x_1,\ldots,x_{10})= \sum_{i=1}^{10}x_iM_i(u,v,w) \in k[u,v,w,x_1,\ldots,x_{10}] , \] 
where $M_1(u,v,w),\ldots,M_{10}(u,v,w)$ are monomials in $u,$ $v$ and $w$ of degree 
$3$ in an arbitrary ordering. 
Stanley's example is given as $A=Q/{\rm Ann}_QF$ corresponding to the polynomial $F$ 
defined above. The algebra $A$ has the Hilbert function 
$(1,13,12,13,1),$ so it does not have 
the strong or weak Lefschetz property. 
More generally, it is shown in \cite{BI} and \cite{Ia} 
that there exist Artinian Gorenstein algebras 
$A$ with a non-unimodal Hilbert function for $\dim A_1 \geq 5.$ 
In Section 5, we will construct Artinian Gorenstein algebras with a unimodal Hilbert 
function which do not have the strong Lefschetz property. 
\end{ex}
\section{Characterization of Lefschetz elements}
In this section we discuss the set of the Lefschetz elements for graded 
Artinian Gorenstein rings $A=k[X_1,\ldots, X_n]/ {\rm Ann}_Q F$ with the standard 
grading. 
\begin{dfn}
Let $G$ be a polynomial in $k[x_1,\ldots,x_n].$
When a family ${\bf B}_d=\{ \alpha^{(d)}_i \}_i$ of homogeneous polynomials of degree 
$d>0$ is given, 
we call the polynomial 
\[ \det \Big( (\alpha^{(d)}_i(X) \alpha^{(d)}_j(X) G(x))_{i,j=1}^{\dim A_d} \Big) 
\in k[x_1,\ldots,x_n] \] 
the $d$-th Hessian of $G$ with respect to ${\bf B}_d,$ 
and denote it by ${\rm Hess}_{{\bf B}_d}^{(d)}G.$ 
We denote the $d$-th Hessian simply by ${\rm Hess}^{(d)}G$ if the choice of 
${\bf B}_d$ is clear. 
\end{dfn}
When $d=1$ and $\alpha^{(1)}_j(X)=X_j,$ $j=1,\ldots ,n,$ the first Hessian 
${\rm Hess}^{(1)}G$ coincides with the usual Hessian: 
\[ {\rm Hess}^{(1)}G= {\rm Hess}\; G:=
\det \left( \frac{\partial^2 G}{\partial x_i \partial x_j}  
\right)_{ij} . \] 

Let us consider the case $A=Q/{\rm Ann}_Q F$ and the higher Hessians of $F$ 
with respect to a $k$-linear basis ${\bf B}_d=\{ \alpha^{(d)}_i \}$ of $A_d.$ 
Note that if we change the $k$-linear basis of $A_d,$ the corresponding 
higher Hessians ${\rm Hess}_{{\bf B}_d}^{(d)}$ 
are just multiplied by nonzero scalars in $k^{\times}.$ 
\begin{thm} {\rm (\cite[Theorem 4]{W2})}
Fix an arbitrary $k$-linear basis ${\bf B}_d$ of $A_d$ for $d=1,\ldots,[D/2].$ 
An element $L=a_1X_1+\cdots +a_nX_n \in A_1$ is a strong Lefschetz element 
of $A=Q/{\rm Ann}_Q F$ 
if and only if $F(a_1,\ldots,a_n) \not= 0$ and 
\[ ({\rm Hess}_{{\bf B}_d}^{(d)}F)(a_1,\ldots,a_n) \not= 0 \]
for $d=1,\ldots,[D/2].$ 
\end{thm}
\begin{proof} Define the identification $[ \;\; ]:A_D \; \tilde{\rightarrow} \; k$ 
by $[\omega(X)] := \omega(X)F(x)$ for any $\omega(X) \in A_D.$ 
Note that $\omega(X)F(x)\in k,$ 
because $\deg \omega = \deg F=D.$ 
Since 
$A$ is a Poincar\'e duality algebra, the necessary and sufficient condition 
for $L=a_1X_1+\cdots +a_nX_n \in A_1$ to be a strong Lefschetz element is that the 
bilinear pairing 
\[ \begin{array}{ccccc}
A_d \times A_d & \rightarrow & A_D & \cong & k \\ 
(\xi,\eta) & \mapsto & L^{D-2d} \xi \eta & \mapsto & [L^{D-2d}\xi \eta ] 
\end{array} \]
is non-degenerate for $d=0,\ldots,[D/2].$  
Therefore $L$ is a Lefschetz element if and only if 
the matrix 
\[ (L^{D-2d}\alpha^{(d)}_i(X)\alpha^{(d)}_j(X)F(x))_{ij} \] 
has nonzero 
determinant. For a homogeneous polynomial $G(x_1,\ldots,x_n) \in k[x_1,\ldots,x_n]$ 
of degree $d,$ we have the formula 
\[ (a_1X_1+\cdots +a_nX_n)^dG(x_1,\ldots,x_n)=d!G(a_1,\ldots,a_n) ,\] 
so 
\[ L^{D-2d}\alpha^{(d)}_i(X)\alpha^{(d)}_j(X)F(x)= (D-2d)!
\alpha^{(d)}_i(X) \alpha^{(d)}_j(X) F(x)|_{(x_1,\ldots,x_n)=(a_1,\ldots,a_n)}. \] 
\end{proof} 
\begin{cor} $(1)$ The algebra $A=Q/{\rm Ann}_Q F$ has the strong Lefschetz 
property if and only if all the higher Hessians 
${\rm Hess}_{{\bf B}_d}^{(d)}F$ with respect to 
a $k$-linear basis ${{\bf B}_d}$ of $A_d,$ $d=1,\ldots,[D/2],$ 
are nonzero polynomials. \\ 
$(2)$ Assume that the socle degree of $A$ is less than $5.$ 
An element $L=a_1X_1+\cdots +a_nX_n$ is 
a strong Lefschetz element if and only if 
\[ F(a_1,\ldots,a_n)\not=0 \quad \textrm{and} \quad 
{\rm Hess}\; F(a_1,\ldots,a_n)\not= 0. \] 
Here ${\rm Hess}\; F$ is the first Hessian of F with respect to a linear 
basis of $A_1.$ 
\end{cor}

\section{Set of Lefschetz elements}
In this section we discuss the set of Lefschetz elements for 
some simple examples of Gorenstein algebras with the strong Lefschetz property 
based on Corollary 3.1. 
\begin{ex} 
Let us consider the Gorenstein ring $A=k[X_1,\ldots, X_n]/ {\rm Ann}_Q F$ 
associated to the Fermat type polynomial 
\[ F=\sum_{i=1}^nx_i^n-n(n-1)s\prod_{i=1}^nx_i, \] 
where $s\in k$ is a parameter. 
One can check that $A$ has the strong Lefschetz property 
for any $s\in k$ as follows. 
For $s=0$ it is easy to see that $A$ has the strong Lefschetz property. 
When $s\not= 0,$ the monomials $\alpha_1:=x_1^d,\ldots,\alpha_n:=x_n^d$ and 
\[ \alpha_I:=\prod_{i\in I}x_i, \quad I\subset \{ 1,\ldots, n \}, \; \# I=d \] 
form a linear basis ${\bf B}_d$ of $A_d$ for $d>1.$ The matrix 
$M^{(d)}=\left( \alpha \beta F \right)_{\alpha,\beta \in {\bf B}_d}$ is of form 
\[ M^{(d)}= \left( \begin{array}{cc} 
M_1^{(d)} & 0 \\ 
0 & M_2^{(d)} 
\end{array} \right), \] 
where $M_1^{(d)}$ is a diagonal matrix of size $n$ with $\det M_1^{(d)}\not= 0,$ 
and $M_2^{(d)}$ is a matrix of 
size ${n \choose d}.$ 
Let us consider the monomial $G:=x_1\cdots x_n$ and the corresponding 
algebra $A':=Q/{\rm Ann}_Q G.$ Then we have 
\[ A' \cong k[X_1,\ldots, X_n]/(X_1^2,\ldots,X_n^2) , \] 
so $A'$ has the strong Lefschetz property. Thus the $d$-th Hessian 
${\rm Hess}^{(d)}G$ 
with respect to the linear basis $\{ \alpha_I \; | \; \# I=d \}$ is nonzero. 
Since 
\[ \det M_2^{(d)}= (-n(n-1)s)^{{n \choose d}}\cdot {\rm Hess}^{(d)}G \not= 0, \] 
we have ${\rm Hess}^{(d)}_{{\bf B}_d}F\not= 0.$ Hence $A$ has the strong 
Lefschetz property. 

We give the explicit 
condition for the Lefschetz element for $n=3,4.$ 
For $n=3$ and $F=x^3+y^3+z^3-6s\cdot xyz,$ $A$ has the following structure:
\[ \begin{array}{ll}
\textrm{Case $s\not=0, 1,$} & A \cong k[X,Y,Z]/(sX^2+YZ,sY^2+XZ,sZ^2+XY), \\ 
\textrm{Case $s=0,$} & A \cong k[X,Y,Z]/(X^3-Y^3,X^3-Z^3,XY,YZ,XZ), \\ 
\textrm{Case $s=1,$} & A \cong k[X,Y,Z]/(X^2+YZ,Y^2+XZ,Z^2+XY,XZ^2,YZ^2). 
\end{array} \] 
Note that $A$ is a complete intersection for $s\not=0,1.$ 
The Hilbert function of $A$ is ${\rm Hilb}(A)=(1,3,3,1)$ for all $s\in k.$ 
The condition for $L=aX+bY+cZ \in A_1$ to be a strong Lefschetz element is that 
\[ a^3+b^3+c^3-6s\cdot abc\not=0 \] 
and 
\[ s^2a^3+s^2b^3+s^2c^3-(1-2s^3)abc\not=0. \]  
It is remarkable that the above condition is described by a single condition 
\[ a^3+b^3+c^3-6s\cdot abc\not=0 \] 
for $s=1/2,(-1\pm\sqrt{-3})/4.$ This is exactly when 
$F$ decomposes into the product of three linear forms. 

For $n=4$ and $F=x^4+y^4+z^4+w^4-12s\cdot xyzw,$ 
the Hilbert function of $A$ is as follows: 
\begin{eqnarray*}
{\rm Hilb}(A)=(1,4,4,4,1), & \textrm{for $s=0,$} \\ 
{\rm Hilb}(A)=(1,4,10,4,1), & \textrm{for $s\not= 0.$} 
\end{eqnarray*}
In this case the condition for 
$L=aX+bY+cZ+dW \in A_1$ to be a strong Lefschetz element is that 
\[ a^4+b^4+c^4+d^4-12s\cdot abcd\not=0 \] 
and 
\[ (1-2s^3-s^4)a^2b^2c^2d^2-2s^3\cdot {\rm symm}(a^5bcd)-s^2\cdot {\rm symm}(a^4b^4) \not=0, \]
where ${\rm symm}(\; \cdot \;)$ means the symmetrization of the indicated monomial.  
\end{ex} 
\begin{ex}
Stanley \cite{S1} studied the strong Lefschetz property of the coinvariant 
algebra of finite Coxeter groups to show the Sperner property for 
the Bruhat ordering on finite Coxeter groups. 
In \cite{MNW}, the set of the Lefschetz elements for the coinvariant 
algebra of the finite Coxeter group is determined except for type $H_4.$ 
Let $V$ be the standard reflection representation of the finite irreducible Coxeter 
group $W.$ Then $W$ acts on the polynomial ring $R={\rm Sym}_{\bf R}V^*$ and 
the $W$-invariant subalgebra $R^W$ is generated by the fundamental 
$W$-invariants $f_1,\ldots,f_r,$ $r=\dim V.$ The coinvariant algebra $R_W$ 
is defined as the quotient algebra $R/(f_1,\ldots,f_r).$ 
It is known that $R_W$ is Gorenstein (see e.g. \cite[Theorem 7.5.1]{Sm}). 
When $W$ is crystallographic, 
$R_W$ is isomorphic to the cohomology ring of the corresponding flag variety. 
In \cite{MNW}, it was shown that the set of Lefschetz elements in $V^*=(R_W)_1$ 
is the complement of the union of the reflection hyperplanes. 
For crystallographic case, their argument is based on the ampleness criterion 
for the ${\bf R}$-divisors on the flag variety, so it is applicable only when 
the field $k$ of coefficients is the field ${\bf R}$ of real numbers. 

Let us consider the case $W=S_3$ and 
\[ R_W={\bf R}[X,Y,Z]/(X+Y+Z,XY+YZ+ZX,XYZ). \]  
The algebra $R_W$ is also given by $R_W={\bf R}[X,Y,Z]/{\rm Ann}\; \Delta$ 
with $\Delta=(x-y)(x-z)(y-z).$ 
The degree one part $(R_W)_1$ has a linear basis ${\bf B}_1=\{ X,Y \}.$ 
Then we have 
\[ {\rm Hess}^{(1)}_{{\bf B}_1} \Delta = -4(x^2+y^2+z^2-xy-yz-zx), \] 
which is a negative definite quadratic form. Hence the set of the Lefschetz 
elements is given by 
\[ \{ (x,y,z) \; | \; \Delta(x,y,z) \not= 0 \} \subset V^*. \] 
If we work in $V^*_{\bf C},$ 
we have to take care of the condition $x^2+y^2+z^2-xy-yz-zx\not= 0,$ too. 
\end{ex}
\section{Gorenstein algebras which do not have the strong Lefschetz property}
The result in Section 3 shows that a polynomial $F$ gives an example of 
Gorenstein algebra which does not have the strong Lefschetz property if 
one of the higher Hessians of $F$ is identically zero. In \cite[Section 4]{HMNW} and 
\cite{W2}, some examples of $F$ with the zero Hessian are discussed. 
In \cite{GN}, it is proved , among other things, that the Hessian of a polynomial in 
$4$ variables does not vanish, unless a variable can be eliminated by means of 
a linear transformation of the variables. The polynomial 
$F=x_0u^2+x_1uv+x_3v^2$ is the simplest example whose Hessian vanishes, but no variables 
can be eliminated by a linear transformation of the variables (see \cite[Example 1]{W2}). 

Here we give examples of forms $F$ such that ${\rm Hess}\; F\not=0$ and ${\rm Hess}^{(2)}F=0.$ 
By using these forms we can also give examples of Gorenstein algebras
$A=Q/{\rm Ann}_Q F$ which do not satisfy the strong Lefschetz property. 
\begin{ex} 
Let us consider the polynomial 
\[ F:= \sum_{j=0}^n x_j^2 u^{n-j}v^j \in k[u,v,x_0,\ldots,x_n] \] 
and the corresponding algebra $A=Q/{\rm Ann}_QF,$ 
where $Q=k[U,V,X_0,\ldots,X_n],$ $U =\partial /\partial u,$ $V =\partial /\partial v$ and 
$X_i =\partial /\partial x_i.$ 
Linear bases of $A_1$ and $A_2$ are given in the following table. 
\[ \begin{array}{|c|c|} 
\hline
 & \textrm{linear basis} \\ 
\hline 
\dim A_1=n+3 & U,V,X_0,\ldots,X_n \\ 
\hline 
\dim A_2= 3n+4 & \alpha_1:=U^2,\alpha_2:=UV,\alpha_3:=V^2,
\alpha_4:=UX_0,\ldots,\alpha_{n+3}:=UX_{n-1}, \\ 
 & \alpha_{n+4}:=VX_1,\ldots,
\alpha_{2n+3}:=VX_n, \\ 
 & \alpha_{2n+4}:=X_0^2,\ldots,\alpha_{3n+4}:=X_n^2 \\ 
\hline
\end{array} \]
It is easy to see that 
$\dim A_0=\dim A_{n+2}=1,$ $\dim A_1= \dim A_{n+1}=n+3$ and 
$\dim A_d=(d+1)n-d^2+2d+4$ for $2\leq d \leq n.$ 
The Hessian of $F$ with respect to the basis above is expressed as follows: 
\begin{align*}  
\lefteqn{{\rm Hess}\; F =} \\ 
 &  2^{n+1}(uv)^{\frac{n(n-1)}{2}} \left\{  
(\sum_{j=0}^{n-1}(n-j)(n-j+1)x_j^2u^{n-j-1}v^j)(\sum_{j=1}^nj(j+1)x_j^2u^{n-j}v^{j-1}) 
\right. \\ &  -uv
\left. (\sum_{j=1}^{n-1}j(n-j)x_j^2u^{n-j-1}v^{j-1})^2 
\right\} \not=0. 
\end{align*} 
On the other hand, we have 
\begin{align*} 
X_i^2U^2F & = 2(n-i)(n-i-1)u^{n-i-2}v^i, \\ 
X_i^2UVF & = 2(n-i)iu^{n-i-1}v^{i-1}, \\ 
X_i^2 V^2F & =  2i(i-1)u^{n-i}v^{i-2}, 
\end{align*}  
and $X_i^2\alpha_j(X)F=0$ for $i=0,\ldots,n$ and $j\geq 4,$ 
so the vectors 
\[ \vec{\xi}_i:=(X_i^2 \alpha_1(X)F,X_i^2\alpha_2(X)F,\ldots,
X_i^2\alpha_{3n+4}(X)F), \quad i=0,\ldots,n, \] 
are linearly dependent. 
Hence we can see that the second Hessian is identically zero, i.e. 
${\rm Hess}^{(2)}F=0$ in $k[u,v,x_0,\ldots,x_n].$ 
This means that the algebra $A$ does not have the strong Lefschetz property. 

For $n=3,$ the algebra $A=Q/{\rm Ann}_Q F$ has the Hilbert function 
$(1,6,13,13,6,1).$  
Since the multiplication map $\times L:A_2 \rightarrow A_3$ cannot be bijective 
for all $L\in A_1,$ $A$ does not have the weak Lefschetz property either. 
\end{ex}
\begin{ex} 
There exists an example of a polynomial $F$ of degree $5$ with $5$ variables such that 
${\rm Hess}\; F\not= 0$ and ${\rm Hess}^{(2)}\; F=0.$ 
Let us choose 
\[ F=x^2u^3+xyu^2v + y^2 uv^2 + z^2v^3 \in k[u,v,x,y,z]. \] 
Then 
\[ U,V,X,Y,Z \in A= k[U,V,X,Y,Z]/{\rm Ann}_Q F \] 
are linearly independent. So we have 
\begin{align*} 
{\rm Hess}\; F = & \; 48u^3v^3 \left( u^5x^4+8u^4vx^3y+16u^3v^2x^2y^2+19u^2v^3x^2z^2 
\right. \\ 
 & \left. +9u^2v^3xy^3+13uv^4xyz^2+2uv^4y^4+4v^5y^2z^2 \right) \not=0. 
\end{align*} 
The monomials 
\[ \alpha_1=U^2,\alpha_2=V^2,\alpha_3=UV,\alpha_4=X^2,\alpha_5=Y^2,\alpha_6=Z^2, \] 
\[ \alpha_7=XY,\alpha_8=UX,\alpha_9=UY,\alpha_{10}=VX,\alpha_{11}=VY,\alpha_{12}=VZ \] 
form a linear basis of $A_2.$ We have 
\[ (\alpha_4\alpha_1F,\alpha_4\alpha_2F,\alpha_4\alpha_3F)= (12u,0,0), \] 
\[ (\alpha_5\alpha_1F,\alpha_5\alpha_2F,\alpha_5\alpha_3F)= (0,4v,4u), \]
\[ (\alpha_6\alpha_1F,\alpha_6\alpha_2F,\alpha_6\alpha_3F)= (0,0,12v), \]
\[ (\alpha_7\alpha_1F,\alpha_7\alpha_2F,\alpha_7\alpha_3F)= (2v,2u,0) \] 
and $\alpha_i\alpha_jF=0$ for $i=4,5,6,7$ and $j\geq 4.$ Hence the vectors 
\[ \vec{\xi}_i=(\alpha_i\alpha_1F,\ldots,\alpha_i\alpha_{12}F), \quad j=4,5,6,7, \] 
are linearly dependent and ${\rm Hess}^{(2)}F=0.$ 
The algebra $A=Q/{\rm Ann}_QF$ has the Hilbert function $(1,5,12,12,5,1).$ 
Since we have ${\rm Hess}^{(2)}F=0,$ the algebra $A$ does not have 
the weak Lefschetz property. 
\end{ex}
\begin{ex} 
The following example is due to Ikeda \cite{Ik}. Let us choose the polynomial 
$F=w^3xy+wx^3z+y^3z^2.$ Then the corresponding algebra $A=Q/{\rm Ann}_QF$ has the 
Hilbert function $(1,4,10,10,4,1).$ We choose the linear bases $X,Y,Z,W$ of $A_1$ 
and 
\[ \alpha_1=W^2,\alpha_2=X^2,\alpha_3=Y^2,\alpha_4=Z^2,\alpha_5=WX, \] 
\[ \alpha_6=WY,\alpha_7=WZ,\alpha_8=XY,\alpha_9=XZ,\alpha_{10}=YZ \] 
of $A_2.$ 
The Hessian is given as follows: 
\begin{align*} 
{\rm Hess}\; F = & \; 8(3w^7xy^4+8w^6x^6-27w^5x^3y^3z+27w^4y^6z^2 \\ 
 &  -45w^3x^5y^2z^2-54w^2x^2y^5z^3+9wx^7yz^3+27x^4y^4z^4) . 
\end{align*}
It is easy to check that the four vectors $(\alpha_6\alpha_iF)_i,(\alpha_7\alpha_iF)_i,(\alpha_8\alpha_iF)_i,
(\alpha_9\alpha_iF)_i$ 
are linearly dependent, so the second Hessian ${\rm Hess}^{(2)}F$ is identically zero. 
\end{ex}
\begin{rem} 
In the above examples, we see that $\dim A_1$ takes each value greater than $3.$ 
It is known that if $\dim A_1 =2,$ the Artinian Gorenstein algebra $A$ 
with the standard grading has the strong Lefschetz property \cite[Proposition 4.4]{HMNW}, 
\cite[Theorem 2.9]{Ia}. 
It is still open whether the Artinian Gorenstein algebra with $\dim A_1 =3$ has 
the strong (or weak) Lefschetz property. 
\end{rem}

\end{document}